\documentclass{article}
\usepackage{amsmath}
\usepackage{amssymb}
\usepackage{graphicx}
\usepackage[all]{xy}
\usepackage{amsthm}

\theoremstyle{plain}
\newtheorem{thm}{Theorem}
\newtheorem{lem}{Lemma}
\newtheorem*{thm1}{Theorem}

\theoremstyle{remark}
\newtheorem*{claim}{Claim}
\newtheorem*{remark}{Remark}

\newcommand{\norm}[1]{\langle \langle #1 \rangle \rangle}
\newcommand{\wt}[1]{\widetilde{#1}}
\newcommand{\co}{\colon\thinspace}

\begin{document}

\title{Surface Groups are Frequently Faithful}
\author{Jason DeBlois and Richard P. Kent IV}

\maketitle

\begin{abstract} 
We show the set of faithful representations of a closed orientable 
hyperbolic surface group is
dense in both irreducible components of the $\mathrm{PSL}_2
(\mathbb{K})$ representation variety, where $\mathbb{K}=\mathbb{C}$ 
or $\mathbb{R}$, 
answering a question of W. Goldman. We also prove the existence of  faithful 
representations
into $\mathrm{PU}(2,1)$ with certain nonintegral Toledo invariants.
\end{abstract}

\section{Introduction}

Let $\pi$ be the fundamental group of a closed oriented surface $\Sigma$ 
of genus $g \geq 2$. For 
$\mathbb{K}=\mathbb{C}$ or $\mathbb{R}$, the set 
$\mathrm{Hom}(\pi,\mathrm{PSL}_2(\mathbb{K}))$
is naturally a $\mathbb{K}$--algebraic set, called the 
\textit{ $\mathrm{PSL}_2(\mathbb{K})$
representation variety} of $\pi$. The representation variety inherits a
topology from its ambient affine space; call this the
\textit{classical topology}. We prove the following
\begin{thm}\label{main}
For $\mathbb{K} = \mathbb{C}$ or $\mathbb{R}$, the set of faithful 
representations is dense in $\mathrm{Hom}(\pi,\mathrm{PSL}_2(\mathbb{K}))$ 
equipped with its classical 
topology.
\end{thm}
\noindent This answers a question of W. Goldman \cite{gold1,gold3}.  

Given a representation $\phi\co \pi \rightarrow \mathrm{PSL}_2(\mathbb{K})$, there is an associated flat principal $\mathrm{PSL}_2(\mathbb{K})$--bundle over $\Sigma$, see \cite{milnor}.  The only obstruction to building a cross-section of this bundle is a class $o_2(\phi)$ in $H^2(\pi; \pi_1(\mathrm{PSL}_2(\mathbb{K}))) \cong \pi_1(\mathrm{PSL}_2(\mathbb{K}))$. When $\mathbb{K}=\mathbb{R}$, $o_2(\phi)$ is the \textit{Euler number} $e(\phi)$---so called as it is the Euler number of the associated $\mathbb R \mathrm P^1$--bundle over $\Sigma$---and when $\mathbb{K}=\mathbb{C}$, $o_2(\phi)$ is the \textit{second Stiefel-Whitney class} $w_2(\phi)$---as it is the second Stiefel-Whitney class of the associated $\mathbb H^3$--bundle over $\Sigma$. The second Stiefel-Whitney class of a real representation is the Euler number modulo two.

In \cite{gold3}, Goldman proves that the topological components of 
$\mathrm{Hom}(\pi, \mathrm{PSL}_2(\mathbb{C}))$ 
are $w_2^{-1}(i)$ for $i=0,1$
and those of $\mathrm{Hom}(\pi, \mathrm{PSL}_2(\mathbb{R}))$ 
are $e^{-1}(n)$ for $2-2g \leq n \leq 2g-2$.
Furthermore, each representation variety has two irreducible components,
corresponding to the two possible Stiefel-Whitney classes \cite{gold1,gold3}.

We prove that in each of these irreducible components, the set of faithful 
representations is a countable intersection of open dense sets, which by 
the Baire Category Theorem is dense. To do this, we use the fact that 
Fuchsian representations have even Euler numbers and the following
\begin{thm}\label{tool} For any nontrivial $w \in \pi$,
there is a representation 
$\phi\co \pi \rightarrow \mathrm{PSL}_2(\mathbb{R})$
such that $e(\phi)=3-2g$ and $\phi(w)$ is nontrivial.
\end{thm}
 
In section 4 we discuss the effect upon the obstruction classes when 
a representation is altered by choosing a new embedding of the entry field 
of its image.  In the final section we use our techniques to prove a
result about representations into $\mathrm{PU}(2,1)$.

\section{Proof of Theorem 2}

\begin{proof}
Fix a hyperbolic structure on $\Sigma$.  If $\alpha$ is a geodesic in $\Sigma$, let $\ell(\alpha)$ denote its length with respect to the chosen hyperbolic structure. 

Let $\omega$ be a closed geodesic in $\Sigma$. Since the closed geodesics
are in 1-1 correspondence with the
conjugacy classes in $\pi_1(\Sigma,p)$ for any choice of basepoint $p$, 
it is enough to show that for such a geodesic, there is 
some basepoint $p$ and a representation of $\pi_1(\Sigma,p)$ sending
some element $w$ in the conjugacy class corresponding to $\omega$ to
a nontrivial isometry of $\mathbb{H}^2$.

Choose a simple closed geodesic $\gamma$ on $\Sigma$ with 
$\ell(\gamma)>2\ell(\omega)$, $\gamma \pitchfork \omega \neq \emptyset$, and 
so that $\gamma$ cuts $\Sigma$ into a surface $\Sigma_a$ of genus $g-1$ 
and a punctured torus $\Sigma_b$.  This can be done by applying high
powers of a pseudo-Anosov homeomorphism to any simple closed curve cutting
off a punctured torus. Let $p \in \gamma$, and let 
$A=\pi_1 (\Sigma_a,p)$ and 
$B=\pi_1 (\Sigma_b,p)$.  Then $\pi_1 (\Sigma,p) = A *_C B$, 
where $C$ is the subgroup generated by $[\gamma]$, the homotopy class of
$\gamma$ relative to $p$.  Given an element $h$ in $\pi_1(\Sigma,p)$, let $\norm{h}$ denote its normal closure in $\pi_1(\Sigma,p)$.

\begin{lem}\label{noncomm} There exists 
$w \in \pi_1 (\Sigma,p)$ in the free homotopy class of $\omega$ such that
\[ w = a_0 b_0 \cdots a_{n-1} b_{n-1}, \]
where, for each $i$, $b_i \in B-\norm{[\gamma]}$ and $a_i \in A-\norm{[\gamma]^k}$ for all $k >1$ and is not conjugate to a power of $[\gamma]$.
\end{lem}

\begin{proof} Orient $\gamma$, 
and write $\omega$ as a concatenation of oriented geodesic arcs,
\[ \omega = \alpha_0 . \beta_0 \hdots \alpha_{n-1} . \beta_{n-1}, \]
where $\alpha_i$ is contained in $\Sigma_a$ and $\beta_i$ is contained
in $\Sigma_b$ for all $i$.

For each $i$ let $\epsilon_i$ be the arc of $\gamma$ from $p$ to the 
terminal point of
$\beta_i$ in the positive direction, and let $\delta_i$ be the 
shortest arc of $\gamma$ from 
the terminal point of $\beta_i$ to its initial point.  For any arc
$\alpha$, let $\overline{\alpha}$ denote $\alpha$ with the opposite
orientation.

Now for $0 \leq i \leq n-1$, define \[
a_i=[\epsilon_{i-1} .\alpha_i .\overline{\delta_i} .\overline{\epsilon_i}]
\ \ \mbox{and}\ \ 
b_i=[\epsilon_i .\delta_i .\beta_i .\overline{\epsilon_i}], \]
where all indices are taken modulo $n$. Then transparently, the word
$w = a_0 b_0 \cdots a_{n-1} b_{n-1}$
is in the free homotopy class of $\omega$.

Note that $a_i$ is freely homotopic to the closed curve
$\alpha_i.\overline{\delta_i}.\overline{\epsilon_i}.\epsilon_{i-1}$ with
length less than $\ell(\omega)+\frac{1}{2}\ell(\gamma)+\ell (\gamma)<
2 \ell (\gamma)$.  This curve is essential and not homotopic into $\gamma$, 
since transverse geodesics intersect 
minimally, and so $a_i$ is not conjugate to a power of $[\gamma]$.
Similarly, note that $b_i$ is freely homotopic to
$\delta_i .\beta_i$. This curve has length less than
$\frac{1}{2}\ell(\gamma)+\ell (\omega)< \ell(\gamma)$ and is essential for 
the same reason as above. 

The lemma follows immediately from the following (compare the proofs of
Lemma 1 and Theorem 6 of \cite{long}):
\begin{claim} If $x \in \norm{[\gamma]^k}$,
then $\ell{(x)} \geq |k| \ell{(\gamma)}$. \end{claim}

By Theorem 1 of \cite{mask}, the cover $\wt{\Sigma}$
corresponding to $\norm{[\gamma]^k}$
is planar.  Let $x$ have the shortest geodesic 
representative among all nontrivial elements of $\norm{[\gamma]^k}$.  Then every 
lift of $x$ to 
$\wt{\Sigma}$ has a simple geodesic representative, for otherwise one 
could find a shorter element in $\pi_1 (\wt{\Sigma}) = \norm{[\gamma]^k}$.  
Suppose
the geodesic representatives of two lifts $\wt{x}_1$ and $\wt{x}_2$ intersect
transversely.  By the planarity of $\wt{\Sigma}$,  the number of points of intersection is even. Choose two such points, and let $\alpha$ be the
shortest of the arcs of $\wt{x}_1$ and $\beta$ the shortest of the arcs
of $\wt{x}_2$ in the complement of these two points.  Then 
$\alpha.\beta$ has a shorter geodesic representative than $\wt{x}_1$, a
contradiction.  This shows that $x$ is some power of a simple closed curve in $\Sigma$.
But by Corollary 4.2 of \cite{hempel}, the only such powers in $\norm{[\gamma]^k}$
are of the form $[\gamma]^{ik}$, $|i| \geq 1$.
\end{proof}

We now define a family of representations $\phi_t\co A *_C B \rightarrow 
\mathrm{PSL}_2 (\mathbb{R})$, $t \in \mathbb{R}$, with Euler number $3-2g$,
using a construction of Goldman \cite{gold3}.
Let
\[ \begin{bmatrix} a & b \\ c & d \end{bmatrix} \in 
   \mathrm{PSL}_2(\mathbb{R}) \]
denote the projective class of 
\[ \begin{pmatrix} a & b \\ c & d \end{pmatrix} \in 
   \mathrm{SL}_2(\mathbb{R}). \] 

Pick two nonparallel simple closed curves on $\Sigma_b$ intersecting 
in a single point, and
let $x$ and $y$ represent them in $B$. Then $B=\langle x,y \rangle$, 
and we define a solvable
representation $\phi_B\co B \rightarrow \mathrm{PSL}_2 (\mathbb R)$ 
by \begin{align*}
  & \phi_B (x) = \begin{bmatrix} \alpha&0 \\ 0&\alpha^{-1} \end{bmatrix} &
  & \phi_B (y) = \begin{bmatrix} \beta&1 \\ 0 &\beta^{-1} \end{bmatrix},
\end{align*}
where $\alpha \neq \pm 1$ and 
$\beta \notin \{\alpha^{r}\,|\,r \in \mathbb{Q}\}$.
We record for future reference that
\[ \phi_B([\gamma]) = \phi_B([x,y])
   = \begin{bmatrix} 1 & \beta(\alpha^2-1) \\ 0 & 1 \end{bmatrix}, \]
is parabolic. 

Choose a Fuchsian representation 
$\phi_A\co A \rightarrow \mathrm{PSL}_2 (\mathbb{R})$ so that 
$\phi_A([\gamma])=\phi_B([\gamma])$. 
Note that since $\phi_A$ is discrete and faithful, and $\phi_A([\gamma])$ 
fixes $\infty$, for no element $a \in A-C$ can $\phi_A(a)$ fix $\infty$.
Therefore, for every $a \in A-C$, $\phi_A(a)$ has a nonzero $2,1$ entry.

To complete the construction, for $t \in \mathbb{R}$ define 
$\lambda_t = \left[ \begin{smallmatrix} 1&t \\ 0&1 \end{smallmatrix} \right]$,
and note that since $\lambda_t$ commutes with elements of $\phi_B(C)$, the
representations $\phi_A$ and $\lambda_t \phi_B \lambda_t^{-1}$ agree
on $C$. The universal property of free products with amalgamation
yields a representation
$\phi_t\co \pi \rightarrow \mathrm{PSL}_2(\mathbb{R})$
given by
\[ \phi_t|_A = \phi_A\ \ \mbox{and}\ \ \phi_t|_B=\lambda_t \phi_B 
  \lambda_t^{-1}. \]
The Euler number of each $\phi_t$ is $3-2g$, just as in the proof of 
Lemma 8.2 of \cite{gold3}.

The following lemma, which uses techniques of the proof of Proposition 1.3 of 
$\cite{shal}$, gives a criterion for an element to
survive $\phi_t$.

\begin{lem}\label{trans}
Suppose $w \in A *_C B$ is of the form
\[ w = a_1 b_1 a_2 b_2 \cdots a_l b_l,\ a_i \in A,\ b_i \in B
    \ \mbox{for}\ 1 \leq i \leq l, \]
where for each $i$, $\phi_0(a_i)$ has a nonzero $2,1$ entry and 
$\phi_0(b_i)$ is hyperbolic. If $t$ is transcendental over the entry field of $\phi_0(\pi)$, then $\phi_t(w)$ is not the identity.
\end{lem}

\begin{proof} Let $t$ be
transcendental over the entry field of $\phi_0(\pi)$. The
lemma follows from the following

\begin{claim} The entries of $\phi_t(w)$
  are polynomials in $t$, where the degree of the $2,2$ entry is $l$,
  the degree of the $1,2$ entry is at most $l$, and the other entries
  have degree at most $l-1$. \end{claim}

We prove the claim by induction on $l$. Suppose $l=1$.  
Let $\left[ \begin{smallmatrix} a&b \\ c&d \end{smallmatrix} \right]=
\phi_A(a_1)$ and $\left[ \begin{smallmatrix}
u & v \\ 0 & u^{-1} \end{smallmatrix} \right]=\phi_B(b_1)$. Then 
\begin{eqnarray*}
 \phi_t(w) & = & \phi_A(a_1)\,\lambda_t \phi_B(b_1) \lambda_t^{-1} \\
         & = & \begin{bmatrix} a & b \\ c & d \end{bmatrix}
               \begin{bmatrix} 1 & t \\ 0 & 1 \end{bmatrix}
               \begin{bmatrix} u & v \\ 0 & u^{-1} \end{bmatrix}
               \begin{bmatrix} 1 & -t \\ 0 & 1 \end{bmatrix} \\
         & = & \begin{bmatrix} a & b \\ c & d \end{bmatrix}
           \begin{bmatrix} u & v+t(u^{-1}-u) \\ 
                 0 & u^{-1} \end{bmatrix} \\
         & = & \begin{bmatrix} au & 
                 av + bu^{-1}+at(u^{-1}-u) \\
                 cu & 
                 cv + du^{-1}+ct(u^{-1}-u) \end{bmatrix}. 
\end{eqnarray*}
Since by assumption $c \neq 0$ and $u \neq \pm 1$ (because 
$\phi_B(b_1)$ is hyperbolic), the claim visibly holds.

Now suppose $l>1$.  Write
\begin{eqnarray*}
 \phi_t(w) &=& \left( \phi_A(a_1)\,\lambda_t \phi_B(b_{1}) 
          \lambda_t^{-1}
          \right) \cdot \left( \phi_A(a_2) \cdots \lambda_t \phi_B(b_l) 
            \lambda_t^{-1} \right) \\
         &=& \begin{bmatrix} p_{1,1}(t) & p_{1,2}(t) \\ p_{2,1}(t) &
               p_{2,2}(t) \end{bmatrix} \begin{bmatrix} q_{1,1}(t) &
               q_{1,2}(t) \\ q_{2,1}(t) & q_{2,2}(t) \end{bmatrix} \\
         &=& \begin{bmatrix} p_{1,1}(t)q_{1,1}(t)+p_{1,2}(t)q_{2,1}(t) &
               p_{1,1}(t)q_{1,2}(t)+p_{1,2}(t)q_{2,2}(t) \\
               p_{2,1}(t)q_{1,1}(t)+p_{2,2}(t)q_{2,1}(t) &
               p_{2,1}(t)q_{1,2}(t)+p_{2,2}(t)q_{2,2}(t) \end{bmatrix}.
\end{eqnarray*}
By the base case, $p_{1,1}$ and $p_{2,1}$ are both constant in $t$, and 
$p_{1,2}$ and $p_{2,2}$ are at most linear in $t$ (where $p_{2,2}$ has a 
nonzero $t$-coefficient).
Then by the inductive hypothesis, $p_{2,2}q_{2,2}$ has degree $l$ and 
$p_{2,1}q_{1,2}$ has degree at most $l-1$, and hence the $2,2$ entry of
$\phi(w)$ has degree $l$ in $t$.  The argument for the other entries is
similar.
\end{proof}

It is easy to see that
the set of elements of $B$ taken by $\phi_B$ to parabolics is precisely
the commutator subgroup of $B$, which is contained in $\norm{[\gamma]}$.  Therefore, by Lemma \ref{noncomm}, there is
a word $w$ in the conjugacy class corresponding to $\omega$ satisfying
the hypotheses of Lemma \ref{trans}.  This proves the theorem.
\end{proof}

\section{Proof of Theorem \ref{main}}

\begin{proof}
Write 
\[ \mathrm{Hom}(\pi, \mathrm{PSL}_2(\mathbb{C}))= X_0 \cup X_1, \]
where $X_i$ is the irreducible component consisting of representations
with Stiefel-Whitney class $i$ for $i=0,1$. Similarly write
\[ \mathrm{Hom}(\pi, \mathrm{PSL}_2(\mathbb{R}))= Y_0 \cup Y_1, \]
where $Y_i$ is the irreducible component consisting of representations
with Euler number equal to $i$ modulo 2 for $i=0,1$.  Note that
$Y_i \subset X_i$ under the natural inclusion
$\mathrm{Hom}(\pi, \mathrm{PSL}_2(\mathbb{R})) \subset
  \mathrm{Hom}(\pi, \mathrm{PSL}_2(\mathbb{C}))$.
For an element $w$ of $\pi$, let $X_w$ (respectively, $Y_w$) be the 
algebraic subset of $X_1$ ($Y_1$) consisting of representations killing $w$.  

By Theorem \ref{tool}, if $w \in \pi$ is nontrivial then 
$X_w \subset X_1$ and $Y_w \subset Y_1$ are proper algebraic subsets.
It is a standard fact about irreducible complex varieties 
(see page 124 of \cite{shaf}) that the complement of any proper
subvariety is dense in the classical topology.  Therefore 
$X_1-X_w$ is an open dense subset of $X_1$.

For arbitrary irreducible real algebraic varieties it is not true that 
the complements of 
proper subvarieties are dense in the classical topology. However it is
true given the additional hypothesis that smooth points are dense (see the 
discussion below), and this follows from Proposition 3.7 of
\cite{gold2}:

\begin{thm1}[Goldman] The set of smooth points of 
$\mathrm{Hom}(\pi,\mathrm{PSL}_2(\mathbb{R}))$ 
is dense in $\mathrm{Hom}(\pi,\mathrm{PSL}_2(\mathbb{R}))$.
\end{thm1}

So $Y_1-Y_w$ is an open dense set of $Y_1$.
But the set of faithful representations in $Y_1$ is precisely
\[ \bigcap_{1 \neq w \in \pi_1 \Sigma} Y_1 - Y_w \]
and similarly for faithful representations in $X_1$.
This is an intersection of open dense subsets, which by the 
Baire Category Theorem is dense.

This proves that the set of faithful representations is dense in the
irreducible component of each representation variety corresponding
to a nonzero Stiefel-Whitney class.  In the other component
this is immediately evident by the above argument, since Fuchsian
representations, for instance, are faithful with Stiefel-Whitney
class 0.
\end{proof}

Here is an example that shows that the complement of a proper subvariety
of an irreducible real algebraic variety is not always dense.
The two-variable polynomial
\[ p(x,y) = y^2 - x^2(x-1) \]
is irreducible, and so the set
\[ V(p) = \{ (x,y) \in \mathbb{R}^2\ |\ p(x,y)=0 \} \]
is an irreducible real algebraic variety.  The points 
$(0,0)$ and $(1,0)$ are the only elements of $V(p)$ with $y$-coordinate
equal to $0$, and if $(x,y) \in V(p)$ with $y \neq 0$, then $y^2 >0$,
which implies that $x>1$.  Hence $(0,0)$ is an isolated point of $V(p)$ in the
classical topology, even though $V(p)$ is one dimensional. It is also
a proper subvariety, which, by the above, is not approached by points
in its complement.

In this example, the origin is not a smooth point of the 
variety. The following fact shows that the situation at smooth points
is analogous to the complex case. We give a proof, adapted from the
proof of the complex case in \cite{shaf}.

\begin{lem}\label{realalg}
Let $X \subset \mathbb{R}^n$ be a real algebraic variety of dimension $k$, 
and let $x \in X$ be a smooth point.  Let $Y \subset X$ be a subvariety 
of dimension $l < k$ with $x \in Y$. Then $x$ is approached by a sequence 
in $X-Y$. \end{lem}

\begin{proof}

The proof is by induction on $k$.

If $k=1$, then $l=0$.  Hence $Y$ consists of a finite collection of
points.  Since $x$ is a smooth point of $X$, there is a chart
$\phi\co (0,1) \rightarrow X$ around $x$.
Now pick any sequence in $(0,1)$ that approaches $\phi^{-1}(x)$
and misses $\phi^{-1}(Y)$.  

Now suppose $k>1$.  For each irreducible component $Y_i$ of $Y$,
Let $y_i \in Y_i-\{x\}$, and choose an affine hyperplane $A$
so that $x \in A$, $y_i \notin A$ for all $i$, and $A$ is transverse to 
$X$ at $x$.
Then $A \cap X$ is $(k-1)$--dimensional and no irreducible component of $Y$ 
is contained in $A \cap X$. So for each $i$,
\[ \mbox{dim}\ Y_i \cap A \leq l-1 < k-1. \]
Furthermore, $x$ is a smooth point of $X \cap A$ since the tangent space at
$x$ is of the proper dimension, and the claim 
follows by induction.
\end{proof}

Now for a variety $X$ in which smooth points are dense, any point $x$ 
in a proper
subvariety $Y$ is approached by smooth points of $X$, each of which
is either in $X-Y$ or approached by a sequence in $X-Y$.  A diagonal
argument gives a sequence in $X-Y$ approaching $x$.

\section{Embeddings of Entry Fields}
Given a representation $\phi\co\pi \rightarrow \mathrm{PSL}_2(\mathbb{\mathbb K})$ with image $\Gamma$, a presentation 
\[
\pi = \langle\,x_1,\hdots ,x_{2g}\,|\,   [x_1,x_2] \cdots [x_{2g-1},x_{2g}] = 1\,\rangle,
\]
and lifts $\wt{\phi}(x_1), \hdots, \wt{\phi}(x_{2g})$ to the universal covering $p\co\wt{\mathrm{PSL}}_2(\mathbb{K}) \rightarrow \mathrm{PSL}_2(\mathbb{K})$, we have
\[
o_2(\phi)=[\wt{\phi}(x_1),\wt{\phi}(x_2)] \cdots [\wt{\phi}(x_{2g-1}),\wt{\phi}(x_{2g})] \in \ker\thinspace p \cong \pi_1 (\mathrm{PSL}_2(\mathbb{K})), 
\]
see \cite{milnor}.

We may choose a new embedding $\sigma$ of the entry field of $\Gamma$ to obtain a new representation $\phi^{\sigma}$ with image $\Gamma^{\sigma}$.  
Computing the Stiefel-Whitney 
class of $\phi^{\sigma}$, we have \begin{eqnarray*}
  w_2(\phi^{\sigma}) & = & 
     [\wt{\phi^{\sigma}}(x_1),\wt{\phi^{\sigma}}(x_2)] 
     \cdots [\wt{\phi^{\sigma}}(x_{2g-1}),\wt{\phi^{\sigma}}(x_{2g})] \\
    & = & (\pm I)^{\sigma} \\
    & = & \pm I, \end{eqnarray*}

\noindent since $\mathbb{Q}$ has a unique 
complex embedding.  We therefore have

\begin{remark} $w_2(\phi^{\sigma})=w_2(\phi)$. \end{remark}

This implies that the Euler numbers of two such real representations have 
the same parity. They
may not be equal, as the following example, from
\cite{macreid} pp.\ 161-162 (originally in \cite{ssw}), shows. 

Consider the Saccheri quadrilateral above, where the angle
at A is $\pi/3$, embedded in the upper half plane so that C is at
$i$ and the side DC lies along the geodesic between $0$ and $\infty$.

\begin{figure}
\begin{center}
\input{quadrilateral.pstex_t}
\end{center}
\end{figure}

The Fuchsian subgroup $F$ of the group generated by reflections in the
sides of this quadrilateral has presentation
\[ \langle\,x,y,z\,|\,x^2=y^2=z^2=(xyz)^3=1\,\rangle, \]
where $x$, $y$, and $z$ are rotations by $\pi$ around B, C, and D, 
respectively. It surjects
$\mathbb{Z}_3 \rtimes (\mathbb{Z}_2 \oplus \mathbb{Z}_2)$
(where $\mathbb{Z}_2 \oplus \mathbb{Z}_2$ acts on $\mathbb{Z}_3$ by
$(1,0) \cdot 1=(0,1) \cdot 1=2$) by
\[ x \mapsto (1,(1,0))\ \ y \mapsto (0,(0,1))\ \ z \mapsto (0,(1,1)). \]
The kernel of this map, $\Gamma_0$, is torsion-free with coarea $4\pi$
and so it is the deck group of a genus 2 surface.  

Let $L$ be the length of the side BC, $c=\cosh^2 L$. As in \cite{macreid}, 
the invariant trace field of $F$ is $\mathbb{Q}(c)$. Note that $\Gamma_0$ is 
contained in the commutator subgroup of $F$, and so its trace field is 
contained in $\mathbb{Q}(c)$.
The element $g_0 = (yx)^2$ is in $\Gamma_0$. 
Conjugate $\Gamma_0$ so that $g=\gamma g_0 \gamma^{-1}$ is
$\left( \begin{smallmatrix} \lambda & 0 \\ 0 & \lambda^{-1}
\end{smallmatrix} \right)$, where $\lambda=2c-1-2\sqrt{c^2-c}$. We may 
conjugate the result to obtain a
group $\Gamma$ whose entries lie in $\mathbb{Q}(c,
\lambda)$---see the discussion on page 115 of \cite{macreid}---and so the
entries of $\Gamma$ lie in $\mathbb{Q}(c,\sqrt{c^2-c})$.

Choose $L$ so that $c=\sqrt{2}$.  Then the entries of $\Gamma$
lie in $\mathbb{Q}(\sqrt{2-\sqrt{2}})$. Let $\sigma$ be the
embedding given by $\sigma(\sqrt{2-\sqrt{2}})=\sqrt{2+\sqrt{2}}$.
Consider the element $h_0=(yz)^2$ of $\Gamma_0$, 
which has trace $1+\frac{c}{c-1}$.
Note that $\sigma(\sqrt{2})=-\sqrt{2}$ and so the trace of $h^{\sigma}$
(where $h=\gamma h_0 \gamma^{-1}$) is $1+\frac{\sqrt{2}}{\sqrt{2}+1}<2$.

Now, by the main theorem of \cite{gold1}, any representation of the
genus 2 surface group with Euler number $\pm 2$ is Fuchsian. Since 
$\Gamma^{\sigma}$ contains the elliptic element 
$h^{\sigma}$, the associated representation must have Euler number 0.

\section{Other Lie Groups}

It is natural to ask to what extent these techniques prove informative
about representations of surface groups into other Lie groups.  One example is the following

\begin{thm}\label{toledo} For $g \geq 3$ and any even $k$, $2<k \leq 2g-2$, there are 
faithful representations $\pi \rightarrow \mathrm{PU}(2,1)$ with Toledo
invariants $\pm (k-4/3)$.
\end{thm}

Theorem 
\ref{toledo} follows from

\begin{thm}\label{blonk} For $g \geq 3$ and any finite
collection of nonidentity elements $\{w_j\} \subset \pi$, there is
a representation $\phi\co \pi \rightarrow \mathrm{PU}(2,1)$ with Toledo
invariant $2g-2-4/3$ such that $\phi(w_j)$ is
nontrivial for all $j$.
\end{thm}

\begin{proof} We recall the construction of Proposition 5.1 of \cite{gkl}. 
Let $\gamma$ be a simple closed geodesic on $\Sigma$ which cuts off
a punctured torus, and let $A =\pi_1(\Sigma_a,p)$ and $B=\pi_1(\Sigma_b,p)$
be the fundamental groups of the punctured genus $g-1$ surface and torus,
respectively, cut off by $\gamma$, where $p$ is a point on $\gamma$. Then 
$A$ and $B$ are free and $\pi_1(\Sigma,p)= A*_C B$, as before.
Choose generators $x_1, \hdots, x_{2g-2}$ for $A$ represented by simple
closed curves on $\Sigma_a$ so that $[x_1,x_2] \cdots [x_{2g-3},x_{2g-2}]=
[\gamma]$ and generators $x_{2g-1}, x_{2g}$ for $B$ represented by simple
closed curves on $\Sigma_b$ so that $[x_{2g-1},x_{2g}]=[\gamma]^{-1}$.  This
determines an isomorphism $\pi \rightarrow A*_C B$.

Now let $Q_a$ be the regular $4(g-1)$-gon in 
$\mathbb{H}^1_{\mathbb{C}} \subset \mathbb{H}^2_{\mathbb{C}}$ with vertex
angles $\pi/(6g-6)$, and let $Q_b$ be the
regular quadrilateral in $\mathbb{H}^2$ with vertex angles $\pi/6$. Define $\phi_a\co A \rightarrow \mathrm{SU}(1,1) \subset
\mathrm{PU}(2,1)$ by sending the $x_i$ to the appropriate transvections 
pairing sides of $Q_a$ so that the quotient of $\mathbb{H}^1_\mathbb{C}$
by $\phi_a(A)$ is a genus $g-1$ surface with one cone point of order 3.
Then $\phi_a([\gamma])$ is an elliptic of order 3 fixing some vertex $v$ of
$Q_a$.  As in \cite{gkl}, this rotates by $2\pi/3$ in 
$\mathbb{H}^1_{\mathbb{C}}$ and by $4\pi/3$ in the normal complex geodesic.
Similarly define $\phi_b\co B \rightarrow \mathrm{SU}(1,1)$ so that the 
quotient is
a punctured torus with one cone point of order $3$, then we may assume without
loss of generality
that $\phi_b([\gamma])=\phi_a([\gamma])$. Let $\theta \in \mathrm{PU}(2,1)$
fix $v$ and interchange $\mathbb{H}^1_{
\mathbb{C}}$ with its normal complex geodesic at that point; then
$\theta \phi_b([\gamma]) \theta^{-1} = \phi_b([\gamma])^{-1}$.  We thus obtain a representation $\phi\co \pi \rightarrow \mathrm{PU}(2,1)$ given by $\phi|_A = \phi_a$ and $\phi|_B = \theta \phi_b \theta^{-1}$.  This
representation has Toledo invariant $2g-2-4/3$ (see \cite{gkl}).

Given a nontrivial element $w \in \pi$, Lemma \ref{noncomm} allows us to choose $\gamma$
so that $w=a_1 b_1 \hdots a_m b_m$ where the $\phi_a(a_i)$ and $\phi_b(b_i)$ are hyperbolic.  In fact, 
it is easy to see that $\gamma$ may be chosen so that the above holds for
any finite collection of elements $\{w_j\}$.

Half the side length of the regular $n$-gon in $\mathbb{H}^1_{\mathbb{C}}$ 
with vertex angles $2\pi/3n$ is $l_n$ given by $\cosh(l_n)=\cos(\frac{\pi}{n})/\sin(\frac{\pi}{3n})$. 
The largest ball around a vertex of the polygon
$Q_a$ which is either disjoint from or identical to each of its
translates under $\phi_a(A)$ has radius $r_a=l_{4(g-1)} \geq l_8$.  
The largest such ball around any vertex of $Q_b$ has radius $r_b=l_4$.  
Let $\Pi_a \co 
\mathbb{H}^2_{\mathbb{C}} \rightarrow \mathbb{H}^1_{\mathbb{C}}$ and $\Pi_b\co \mathbb{H}^2_{\mathbb{C}}
\rightarrow \theta(\mathbb{H}^1_{\mathbb{C}})$ be the orthogonal 
projections. An elementary argument, similar to the proof of Proposition 3.1 of \cite{gkl}, shows that
$S_a=\Pi_a^{-1}(\mathbb{H}^1_{\mathbb{C}}-B(r_a,v))$ 
does not intersect $S_b = 
\Pi_b^{-1}(\theta(\mathbb{H}^1_{\mathbb{C}}-B(r_b,v)))$---this is where the hypothesis that $g \geq 3$ is used, as $\Pi_a^{-1}(\mathbb{H}^1_{\mathbb{C}}-B(l_4,v))$ and $\Pi_b^{-1}(\theta(\mathbb{H}^1_{\mathbb{C}}-B(l_4,v)))$ do indeed intersect.
Therefore any hyperbolic element of $\phi_a(A)$ takes $S_b$ inside $S_a$
and vice-versa. By a ping-pong argument, $\phi(w)$ is 
a nontrivial isometry of $\mathbb{H}^2_{\mathbb{C}}$.
\end{proof}

\begin{proof}[Proof of Theorem \ref{toledo}]

By the main theorem of \cite{xia}, there is a unique topological component $\mathrm T$
of the $\mathrm{PU}(2,1)$ representation variety of $\pi$ consisting of
representations with Toledo invariant $2g-2-4/3$. 
Let $W_1, \hdots, W_n$
be the irreducible components
that contain smooth points in $\mathrm T$.
Suppose that for each $i$ there is an element $w_i \neq 1$ such that all 
representations in $W_i$ kill $w_i$.  Theorem \ref{blonk} produces a 
representation in which all of the $w_i$ survive; furthermore, this
representation is a smooth point of $\mathrm T$ (this follows from \cite{gold2}, 
Proposition 3.7), a contradiction.  Hence there is some irreducible
component $W$ for which each $X_{w}= \{\rho \in W\ |\ \rho(w)=1\}$ is
a proper subvariety, and a representation $\rho \in W$ with 
Toledo invariant $2g-2-4/3$ which is a smooth point of $W$.  Let 
$U$ be an open neighborhood of $\rho$ in $W$ consisting entirely of
smooth points.  Then for each $w \neq 1$, $X_w \cap U$ is nowhere
dense in $U$ by Lemma \ref{realalg}.  The union of such sets is nowhere dense in $U$, and so faithful representations exist in $U$.

The more general statement follows from the fact that, given a surface
$\Sigma$ of genus $g \geq 3$ and a finite set of nontrivial elements in
$\pi_1(\Sigma)$, there is a degree one map from $\Sigma$ to a surface of
genus $g-1$ so that each element of the finite set survives the induced
map between fundamental groups.
\end{proof}

\begin{remark}
A stronger theorem concerning $\mathrm{PU}(2,1)$ has been announced by N. Gusevskii. 
\end{remark}

\section*{Acknowledgements} The authors thank Peter Storm
and Ben McReynolds for making us aware of this question and William
Goldman for providing a copy of \cite{gold1}.
We also thank Alan Reid for useful conversations, in particular for making 
us aware of \cite{shal}.

\end{document}